\documentclass[11pt]{article}
\usepackage{amsmath,amscd,amsthm,amsfonts,amssymb, latexsym}
\usepackage{eucal}

\newtheorem{lemma}{Lemma}

\newtheorem{theorem}{Theorem}

\theoremstyle{definition}

\newcommand{\sltwo}{\ensuremath{\mathfrak{sl_{2}}}}
\newcommand{\g}{\ensuremath{\mathfrak{g}}}
\newcommand{\n}{\ensuremath{\mathfrak{n}}}
\newcommand{\h}{\ensuremath{\mathfrak{h}}}

\renewcommand{\b}{\ensuremath{\mathfrak{b}}}

\newcommand{\Lg}{\ensuremath{L\mathfrak{g}}}
\newcommand{\Lns}{\ensuremath{L \n_{+}^{\sigma}}}
\newcommand{\Lgs}{\ensuremath{L \g^{\sigma}}}

\newcommand{\ghat}{\ensuremath{\hat{\mathfrak{g}}}}
\newcommand{\ghats}{\ensuremath{\ghat^{\sigma}}}

\newcommand{\Vg}{\ensuremath{V_{k}(\g)}}

\newcommand{\M}{\ensuremath{M}}
\newcommand{\Ms}{\ensuremath{\M^{\sigma}}}
\newcommand{\Msbar}{\ensuremath{\widetilde{M}^{\sigma}}}
\newcommand{\fock}{\ensuremath{M^{\sigma} \otimes \pi^{k+\dcox,\sigma}_{\chi}}}
\newcommand{\Yp}{\ensuremath{Y_{\pi}}}
\newcommand{\Yps}{\ensuremath{Y^{\sigma}_{\pi}}}

\newcommand{\Ym}{\ensuremath{Y_{M}}}
\newcommand{\Yms}{\ensuremath{Y^{\sigma}_{M}}}
\newcommand{\Ymp}{\ensuremath{Y_{M \otimes \pi}}}
\newcommand{\Ymps}{\ensuremath{\Ymp^{\sigma}}}
\newcommand{\Yg}{\ensuremath{Y_{\Vg}}}
\newcommand{\Ygs}{\ensuremath{\Yg^{\sigma}}}

\newcommand{\A}{\ensuremath{\mathcal{A}}}
\newcommand{\As}{\ensuremath{\A^{\sigma}}}

\renewcommand{\H}{\ensuremath{\mathcal{H}}}
\newcommand{\Hs}{\ensuremath{\H^{\sigma}}}

\newcommand{\vac}{\ensuremath{\rhd}}

\newcommand{\sigg}{\ensuremath{\sigma_{\g}}}
\newcommand{\sigLg}{\ensuremath{\sigma_{\Lg}}}
\newcommand{\sigM}{\ensuremath{\sigma_{\M}}}
\newcommand{\sigpi}{\ensuremath{\sigma_{\pi}}}
\newcommand{\sigMpi}{\ensuremath{\sigma_{\M \otimes \pi}}}
\newcommand{\sigVg}{\ensuremath{\sigma_{\Vg}}}

\newcommand{\wak}{\ensuremath{W_{\chi, k}}}
\newcommand{\waks}{\ensuremath{W^{\sigma}_{k,\chi}}}
\newcommand{\rwaks}{\ensuremath{\overline{W}^{\sigma}_{\chi}}}
\newcommand{\crwaks}{\ensuremath{\overline{W}^{\sigma,*}_{\chi}}}

\newcommand{\verm}{\ensuremath{M_{k,\chi}}}

\newcommand{\znz}[1]{\ensuremath{\mathbb{Z}/{#1}\mathbb{Z}}}

\newcommand{\dcox}{\ensuremath{\check{h}}}

\begin{document}

\title{Wakimoto Modules for Twisted Affine Lie Algebras}
\author{Matthew  Szczesny \thanks{The author is partially supported by an NSERC graduate fellowship} \\
	\\
	Department of Mathematics \\
	 University of California, Berkeley \\ 
	CA, 94720, USA }
\date{June 8, 2001}
\maketitle

\begin{abstract}
We construct Wakimoto modules for twisted affine Lie algebras, and interpret this 
construction in terms of vertex algebras and their twisted modules. Using the Wakimoto 
construction, we prove the Kac-Kazhdan conjecture on the characters of irreducible modules
with generic critical highest weights in the twisted case.  We provide  explicit formulas for the 
twisted fields in the case of $A_{2}^{(2)}$. 
\end{abstract}

\section{Introduction} \label{S:intro}

Wakimoto modules are an important family of highest weight representations of affine Lie algebras. 
They were first constructed in \cite{WAK} for $\hat{\mathfrak{sl_{2}}}$, and in \cite{FF1} for an arbitrary untwisted 
affine lie algebra $\ghat$. The Wakimoto module $\wak$ of level $k$ and highest weight $\chi$ 
has the same character as the corresponding Verma module $\verm$, yet may possess a 
different composition series. Wakimoto modules have important applications in conformal field theory,
representation theory, integrable systems, as well as other areas. 

In this paper, we construct Wakimoto modules for twisted affine Lie algebras. These should have 
analogous applications, for example in conformal field theory, where they can be used to 
bosonize orbifold models.
We compute explicit realizations of these representations for the twisted affine lie algebra $A_{2}^{(2)}$.

It is our goal to cast the Wakimoto construction in the language of vertex algebras and their modules. 
In the untwisted case, the picture is as follows. To $\ghat$ at level $k$, we can associate the vacuum
module vertex algebra
$\Vg$. The Wakimoto construction maps $\Vg$ to a tensor product of two vertex algebras $M \otimes \pi^{k+\dcox}$ , where
$M$ is a $\beta \gamma$ system, and $\pi^{k+\dcox}$ a Heisenberg vertex algebra. $\pi^{k+\dcox}$ has a family of modules 
$\pi^{k+\dcox}_{\chi}$,
and so $\wak = M \otimes \pi^{k+\dcox}_{\chi}$, an $M \otimes \pi^{k+\dcox}$ module, inherits a structure of $\Vg$ module via the
homomorphism $\iota : \Vg \rightarrow M \otimes \pi^{k+\dcox}$.

In the twisted case, a similar picture holds. Let $\sigg$ be a diagram automorphism of $\g$, and let
$\ghats$ be the corresponding twisted affine Lie algebra. $M$ has a twisted module $M^{\sigma}$, and
$\pi^{k+\dcox}$ has a family of twisted modules $\pi^{k+ \dcox, \sigma}_{\chi}$. Thus we obtain a family of twisted 
$M \otimes \pi^{k+\dcox}$ - modules $\waks=M^{\sigma} \otimes \pi^{k+\dcox,\sigma}_{\chi}$. 
$\waks$, when viewed as a twisted $\Vg$ - module via the embedding $\iota: \Vg \rightarrow M \otimes \pi^{k+\dcox}$,
has the the structure of a $\ghats$ module of level k. 

Once Wakimoto modules for twisted affine algebras have been constructed, they can be used to prove
the Kac-Kazhdan conjecture about characters of irreducible modules at the critical level $k= - \dcox$. In
the final sections of this paper we state and prove this conjecture in the twisted case. Our proof
follows the one in \cite{F1} for untwisted affine algebras.  

I would like to thank Edward Frenkel for suggesting this problem and for many very valuable discussions. 

\section{Twisted Affine Lie Algebras}	\label{twistaff}

In this section we briefly review the definition of twisted affine algebras following \cite{K1}. 
We begin by introducing some standard notation. Given a Kac-Moody algebra $\mathfrak{k}$, we use the 
notation $\Lambda(\mathfrak{k})$, $\Delta(\mathfrak{k})$, $\Delta^{+}(\mathfrak{k})$, 
$\Delta^{+}_{re}(\mathfrak{k})$, and $\Delta^{+}_{im}(\mathfrak{k})$, to denote respectively, the root 
lattice of $\mathfrak{k}$, the set of non-zero roots of $\mathfrak{k}$, the positive roots, 
the positive real roots, and the positive imaginary roots. 

Let $\g$ be a finite-dimensional complex simple Lie algebra, possessing a non-trivial diagram automorphism $\sigg$ 
of finite order $N$. $N=2$ when $\g = A_{n}, D_{n}, E_{6}$ and $N=3$ when $\g = D_{4}$. We will use
the notation $\sigg$ both to refer to the action of the automorphism on $\g$, as well as the action on the root
lattice of $\g$. Let $\epsilon = e^{\frac{2 \pi i}{N}}$, and let us write
\[
	\g = \bigoplus_{j \in \znz{N}} \g_{j}
\]
where $\g_{j}$ is the $\epsilon^{j}$ eigenspace of $\sigma$. Furthermore, if $\mathfrak{l} \subset \g$
is a subspace in $\g$, let $\mathfrak{l}_{i} = \mathfrak{l} \cap \g_{i}$. Let
\[
	\Lg^{\sigma} = \mathbb{C} \cdot d \oplus  \bigoplus_{j =0}^{N-1} \g_{j} \otimes t^{\frac{j}{N}} \mathbb{C}((t))
\]
with Lie bracket 
\begin{align*}
[x \otimes t^{n}, y \otimes t^{m}] &= [x,y] \otimes t^{n+m} \\ 
	[d,x \otimes t^{n}] &= n (x \otimes t^{n})
\end{align*}
Let $\check{r}$ denote the dual Coxeter number of $\g$. 
$\Lgs$ has a central extension $\ghats = \Lgs \oplus \mathbb{C} \cdot K$ defined by the cocyle
\[
	\omega (x \otimes t^{n}, y \otimes t^{m}) = n (x,y) \delta_{n,-m} K
\]
where $(x,y) = \frac{1}{2 \check{r}} Tr (ad(x) ad(y))$ denotes the normalized invariant inner product on $\g$.
In our discussion we will also use the untwisted loop algebra
\[
	\Lg = \mathbb{C} \cdot d \oplus \g \otimes \mathbb{C}((t))
\]
with the same commutation relations as above, and its central extension $\ghat$ with respect to the cocyle $\omega$.
Given $X \in \g$, we will henceforth use the notation $X_{n}$ to denote $X \otimes t^{n}$.
The automorphism $\sigg$ of $\g$ induces an automorphism $\sigLg$ of $\Lg$ and (which extends to $\ghat$), which acts
by $\sigLg(X_{n}) = \sigg(X)_{n}$. We denote the dual Coxeter number of $\ghats$ by $\dcox$. It is equal to the 
dual Coexter number of $\ghat$.

Let $\g = \n_{-} \oplus \h \oplus \n_{+}$ be the triangular decomposition of $\g$. $\sigg$ preserves it, and we have
$\Lgs = L \n^{\sigma}_{-} \oplus L \h^{\sigma} \oplus \Lns$. The Lie subalgebra $\g_{0}$ fixed by $\sigg$ is
simple, and each $\g_{j}$, $(j=0, \cdots, N-1)$ is an irreducible representation of $\g_{0}$. In the case $N=3$, 
the $\g_{0}$ representations $\g_{1}$ and $\g_{-1}$ are equivalent. Let 
$\Delta_{j}$ denote the set of non-zero weights occurring in $\g_{j}$ (viewed as a $\g_{0}$ representation), 
and $\g_{j,\alpha}$ 
$, \, \alpha \in \Delta_{j}$, the corresponding weight space. We have that $\Delta_{j} \subset \Lambda(\g_{0})$,
$dim (\g_{j,\alpha}) = 1$, and $[\g_{i,\alpha},\g_{j,\beta}] \subset g_{(i+j) \, mod \,N , \alpha + \beta}$. 
Let us fix a basis $E_{j,\alpha}$ for $\g_{j,\alpha}, \alpha \in \Delta_{j}$, and $H_{j,a}$ for $\h_{j}$,
$a = 1, \cdots, dim(\h_{j})$. Then $\{ E_{n,(j,\alpha)}, H_{n,(j,a)}, d \}$, $n \in \mathbb{Z}$ 
(resp. $n \in \frac{j}{N} + \mathbb{Z}$) forms a basis for $\Lg$ (resp. $\Lgs$).

Let $\n_{j,\alpha} = \g_{j,\alpha} \cap \n_{+}$,
and let $\Delta^{+}_{j} = \{ \alpha \in \Delta_{j} \vert \n_{j,\alpha} \ne 0 \} $ 
We therefore obtain that
\[
	\Lns = \bigoplus_{j=0}^{N-1} \bigoplus_{\alpha \in \Delta^{+}_{j}} \n_{j,\alpha} \otimes t^{\frac{j}{N}} \mathbb{C}((t))
\]
Thus $E_{n,(j,\alpha)}, \, n \in \frac{j}{N} + \mathbb{Z}, \alpha \in \Delta^{+}_{j}$ forms a basis for
$\Lns$.

\section{Vertex Algebras and Twisted Modules}

For an introduction to vertex algebras see (\cite{FB}, \cite{K2}). For a definition of twisted module see \cite{D} and \cite{BKT}. 
We will denote the vacuum vector by $\vac$, and use $\sim$ to denote the singular terms of the operator product expansion (OPE). 
If $V$ is a vertex algebra, and $A \in V$ is a vector of conformal
weight $p$, then we write $Y(A,z) = \sum_{n \in \mathbb{Z}} A[n] z^{-n-p}$ - i.e. $A[n]$ denotes the Fourier
coefficient in front of $z^{-n-p}$. 
It is was shown in (\cite{FZ}, \cite{FF3})that 
\[
	\Vg = Ind^{\ghat}_{\g[[t]]} \mathbb{C}_{k}
\]
(where $\mathbb{C}_{k}$ denotes the one-dimensional representation on which $\g[[t]]$ acts trivially
and $K$ acts by $k \in \mathbb{C}$) has a vertex algebra structure. Given $J \in \g$, this 
structure assigns the field
\[
	J(z) = \Yg(J_{-1} \vac, z) =  \sum_{n \in \mathbb{Z}} J_{n} z^{-n-1}  
\]
The OPE is as follows:
\begin{equation}
	J^{1}(z) J^{2}(w) \sim \frac{[J^{1},J^{2}](w)}{z-w} + \frac{k(J^{1},J^{2})}{(z-w)^{2}}
\end{equation}
The automorphism $\sigg$ of $\g$ induces a vertex algebra automorphism $\sigVg$ of $\Vg$ which acts by
\[
	\sigVg (J^{1}_{n_{1}} \cdots J^{k}_{n_{k}} \vac) \rightarrow \sigg(J^{1})_{n_{1}} \cdots \sigg(J^{k})_{n_{k}} \vac, \qquad n_{i} < 0
\]
\noindent To $J \in \g_{j}$ we can also assign the twisted field
\[
	J^{\sigma}(z) = \Ygs(J_{-1} \vac, z) = \sum_{n \in \frac{j}{N} + \mathbb{Z}} J_{n} z^{-n-1}
\]
\noindent If $J^{1} \in \g_{i}, J^{2} \in \g$, the OPE of the twisted fields is as follows:
\begin{equation} 
	J^{1,\sigma}(z)J^{2,\sigma}(w) \sim \frac{z^{-i/N}w^{i/N}[J^{1},J^{2}]^{\sigma}(w)}{(z-w)} +
		 k(J^{1},J^{2}) \partial_{w} \frac{z^{-i/N}w^{i/N}}{(z-w)^{2}} \label{twistedOPE}
\end{equation}

The commutation relations between the Fourier coefficients of fields
can be recovered from the singular terms in the OPE.

\section{Twisted Fock Spaces}

We proceed to define various Bosonic vertex algebras and their twisted modules that will be used 
in the twisted Wakimoto realization.

\subsection{The $\beta \gamma$ System}
 
Let $\A$ (respectively $\As$) denote the Heisenberg Lie algebra with generators 
$ \{ \tilde{a}^{*}_{n,(j,\alpha)},$ $ \tilde{a}_{n,(j,\alpha)}, \mathbf{1} \}$ $, n \in \mathbb{Z}, 
j = 0, \cdots, N-1,$ $ \alpha \in \Delta^{+}_{j}$,
(respectively $\{ a^{*}_{n,(j,\alpha)},$ $ a_{m,(j,\alpha)}, \mathbf{1} \}$ $, n \in - \frac{j}{N} + \mathbb{Z}, $ $ m \in \frac{j}{N} + \mathbb{Z}, 
j = 0, \cdots, N-1, \alpha \in \Delta^{+}_{j}$)
and commutation relations
\[
	[\tilde{a}_{n,(j,\alpha)}, \tilde{a}^{*}_{m,(k,\beta)} ] = \delta_{n,-m} \delta_{j,k} \delta_{\alpha, \beta} \mathbf{1}
\] 
(respectively, the same commutation relations with $\tilde{a}, \tilde{a}^{*}$ 's replaced with $a, a^{*}$ 's). The 
element $\mathbf{1}$ is central.
$\A$ (resp. $\As$) has a ``positive'' abelian Lie subalgebra $\A_{+}$ (resp. $\As_{+}$) spanned by 
$\tilde{a}^{*}_{n,(j,\alpha)}, n > 0$, $\tilde{a}_{m,(j,\alpha)}, m \ge 0$ (resp. 
$a^{*}_{n,(j,\alpha)}, n > 0$, $a_{m,(j,\alpha)}, m \ge 0$ )

Let $\mathbf{C}$ denote the 1-dimensional representation of $\A_{+} \oplus \mathbb{C} \cdot \mathbf{1}$ (resp. $\As_{+} \oplus \mathbb{C} \cdot \mathbf{1}$) 
on which $\A_{+}$ (resp. $\As_{+}$) acts trivially, and $\mathbf{1}$ acts by $1$, and let
\[
	\M = Ind_{\A_{+} \oplus \mathbb{C} \cdot \mathbf{1} }^{\A} \mathbf{C}
\]
and 
\[
	\Ms = Ind_{\As_{+} \oplus \mathbb{C} \cdot \mathbf{1}}^{\As} \mathbf{C}
\]

It is well-known that $\M$ has a vertex algebra structure (see \cite{FB}), generated in the sense of the reconstruction theorem
(for the reconstruction theorem for vertex algebras see \cite{K2}, \cite{FKRW})
by the field assignments (here $\vac$ denotes the vacuum vector):
\[
	 \tilde{a}^{*}_{(j,\alpha)}(z) = Y_{\M}(\tilde{a}^{*}_{0,(j,\alpha)} \vac, z) = \sum_{n \in \mathbb{Z}} \tilde{a}^{*}_{n,(j,\alpha)} z^{-n} 
\]
\[
	 \tilde{a}_{(j,\alpha)}(z) = Y_{\M}(\tilde{a}_{-1,(j,\alpha)} \vac, z) = \sum_{n \in \mathbb{Z}} \tilde{a}_{n,(j,\alpha)} z^{-n-1}
\]
$\M$ has a vertex algebra automorphism of finite order $N$, which we will denote $\sigM$, 
induced from the automorphism of $\A$ taking
\begin{align*}
	\tilde{a}^{*}_{n,(j,\alpha)} & \rightarrow \epsilon^{-j} \tilde{a}^{*}_{n,(j,\alpha)} \\
	\tilde{a}_{n,(j,\alpha)} & \rightarrow \epsilon^{j} \tilde{a}_{n,(j,\alpha)} 
\end{align*}

$\Ms$ carries the structure of a $\sigma_{M}$ - twisted $M$-module.
The vertex operation
\[
	Y_{\M}^{\sigma} : \M \rightarrow End(\Ms)[[z^{\frac{1}{N}},z^{-\frac{1}{N}}]]
\]
is generated by the fields
\[
	a_{(j,\alpha)}^{*} (z) =  Y_{\M}^{\sigma}(\tilde{a}^{*}_{0,(j,\alpha)} \vac, z) = \sum_{n \in \frac{-j}{N} + \mathbb{Z}} a^{*}_{n,(j,\alpha)} z^{-n} 
\]
\[
	 a_{(j,\alpha)}(z) = Y_{M}^{\sigma}(\tilde{a}_{-1,(j,\alpha)} \vac, z) = \sum_{n \in \frac{j}{N} + \mathbb{Z}} a_{n,(j,\alpha)} z^{-n-1}
\]

\subsection{The Free Boson} \label{FreeBoson}

Let $\H$ (resp. $\Hs$) denote the Heisenberg Lie algebra with generators $\{ \tilde{b}_{n,(i,a)}, \mathbf{1} \}$, $\, n \in \mathbb{Z},$ $i= 0, \cdots, N-1, 
$ $a = 1, \cdots, dim(\h_{i})$,
(resp.  $ \{ b_{n,(i,a)}, \mathbf{1} \},$ $\,  n \in \frac{i}{N} + \mathbb{Z},$ $ i=0 \cdots N-1,$ $a = 1, \cdots, dim(\h_{i}) $) and commutation relations
\[
	[\tilde{b}_{n, (i,a)}, \tilde{b}_{m, (j,b)}] = n (H_{i,a},H_{j,b}) \delta_{n,-m} \mathbf{1}
\]
(resp. $\tilde{b}$'s replaced by $b$'s). The element $\mathbf{1}$ is central. 
$\H$ (resp. $\Hs$) has an abelian subalgebra $\H_{+}$ (resp. $\Hs_{+}$) spanned by $\tilde{b}_{n,(i,a)}, \,  n \geq 0$ 
(resp. $b_{n,(i,a)}, \, n \geq 0$). Let $\mathbf{C}_{r}$ denote the one-dimensional representation of $\H_{+} \oplus \mathbb{C} \cdot \mathbf{1}$ 
on which $\H_{+}$ acts trivially and  $\mathbf{1}$ acts by $r$, and let
\[
	\pi^{r} = Ind_{\H_{+} \oplus \mathbb{C} \cdot \mathbf{1} }^{\H} \mathbf{C}_{r}
\]
Then it is well-known that $\pi^{r}$ has a vertex algebra structure (see \cite{FB}), generated by the fields
\[
	\tilde{b}_{i,a}(z) = Y_{\pi}(\tilde{b}_{-1,(i,a)} \vac, z) = \sum_{n \in \mathbb{Z}} \tilde{b}_{n,(i,a)} z^{-n-1}
\] 

Let $\chi \in \h_{0}^{*}$, and let $c_{a} = \chi(H_{0,a})$. Let $\mathbf{C}_{r,\chi}$ denote the one-dimensional representation 
of $\Hs_{+}$ where $b_{0,(0,a)}$ acts
by $c_{a}$, all other generators of $\Hs_{+}$ act by $0$, and $\mathbf{1}$ acts by $r$. Let

\[
	\pi^{r,\sigma}_{\chi} = Ind_{\Hs_{+} \oplus \mathbb{C} \cdot \mathbf{1}}^{\Hs} \mathbf{C}_{r,\chi}
\]
$\pi^{r}$ has a vertex algebra automorphism of order $N$, which we will denote $\sigpi$, that is induced
from the automorphism of $\H$ sending $\tilde{b}_{n,(i,a)} \rightarrow \epsilon^{i} \tilde{b}_{n,(i,a)}$. 
Then $\pi^{r,\sigma}_{\chi}$ has a structure of $\sigpi$ - twisted $\pi^{r}$ - module, where the vertex operation
$Y_{\pi}^{\sigma}: \pi^{r} \rightarrow End(\pi^{r,\sigma}_{\chi})[[z^{\frac{1}{N}}, z^{\frac{-1}{N}}]]$ is generated by the fields
\[
	b_{i,a}(z) = Y_{\pi}^{\sigma}(\tilde{b}_{-1, (i,a)} \vac, z) = \sum_{n \in \frac{i}{N} + \mathbb{Z}} b_{n,(i,a)} z^{-n-1}
\]

\subsection{Tensor Products of Fock Spaces}

$M \otimes \pi^{r}$ has a vertex algebra structure with the vertex operation
\[
	\Ymp(A \otimes B, z) = \Ym(A,z) \otimes \Yp(B,z), \qquad A \in M, B \in \pi^{r}
\]
and $\sigMpi=\sigM \otimes \sigpi$ is an automorphism of $M \otimes \pi^{r}$ of order $N$.
$\Ms \otimes \pi^{r,\sigma}_{\chi}$ has the structure of a $\sigMpi$ - twisted $\M \otimes \pi$ - module, with the vertex operation
given by
\[
	\Ymps(A \otimes B, z) = \Yms(A,z) \otimes \Yps(B,z), \qquad A \in M, B \in \pi^{r}
\]

\section{The Twisted Wakimoto Construction}

In this section, we prove the existence of twisted Wakimoto modules using the untwisted construction. By
an embedding of vertex algebras we mean a vertex algebra homomorphism which is both injective and surjective.   
The following theorem follows from (\cite{FF2}, \cite{LEC}) (see also \cite{FB} for the case $\g=\sltwo$)

\begin{theorem} 
(Untwisted Wakimoto Construction) (i) There exists an embedding of vertex algebras \label{untwistwak}
\begin{equation}
\begin{CD}
	\Vg @>{\iota}>> \M \otimes \pi^{k+\dcox}	\label{voaembedding}
\end{CD}
\end{equation}	
for any $k \in \mathbb{C}$ \\
(ii) When $k= - \dcox$, there exists a vertex algebra homomorphism
\begin{equation}
\begin{CD}
	V_{- \dcox}(\g) @>{\kappa}>> \M 	\label{voacritembedding}
\end{CD}
\end{equation}
This homomorphism is neither surjective nor injective. 	
\end{theorem}

Note that when $k = - \dcox$ (critical level), $V_{- \dcox}(\g)$ is not a conformal vertex algebra. \\

\noindent \textbf{Remark:} Given vector spaces $V, W$, each equipped with an action of $\znz{N}$,
we will refer to a linear map $f: V \rightarrow W$ simply as ``equivariant'' if $f$ commutes with
this action.  \\

\noindent  $\znz{N}$ acts on $\Vg$ via $\sigVg$ and on $\M \otimes \pi^{k+\dcox}$ via $\sigMpi$. We can therefore 
ask whether the vertex algebra homomorphisms $\iota$ and $\kappa$ are equivariant. This is indeed the case:

\begin{theorem}
	$\iota$ and $\kappa$ are equivariant vertex algebra homomorphisms.
\end{theorem}
\begin{proof}

Let $G$ be the simply-connected algebraic group with Lie algebra $\g$. $\sigg$ descends to 
$G$, preserving $B_{-}$, the Borel subgroup with Lie algebra $\b_{-}$. Thus $\sigg$
descends to the Flag Manifold $G/B_{-}$, and in particular preserves the big cell $N_{+} \cdot B_{-} \cong \n_{+}$.
This induces an action of $\znz{N}$ on $Vect(N_{+} \cdot B_{-})$ - the Lie algebra of vector fields on $N_{+} \cdot B_{-}$, 
denoted $\sigma_{N_{+}}$.
It follows from this that the embedding of $\g$ into $Vect(N_{+})$ induced by the left 
action of $G$ on $G/B_{-}$ is equivariant.

Let $\{E_{\alpha}, F_{\alpha} \}, \, \alpha \in \Delta^{+}(\g)$ be a basis of $\n_{-} \oplus \n_{+}$ such that
$E_{\alpha} \in \g_{\alpha}, F_{\alpha} \in \g_{-\alpha}$, and $\sigg(E_{\alpha}) = E_{\sigg(\alpha)}, 
\sigg(F_{\alpha}) = F_{\sigg(\alpha)}$. Let $\{ \alpha_{1}, \cdots, \alpha_{r} \}$ be the set of
simple roots, and set $H_{\alpha_{i}} = [E_{\alpha_{i}}, F_{\alpha_{i}}]$. Then $\{E_{\alpha}, F_{\alpha}, H_{\alpha_{i}} \}$   
is a basis for $\g$.

We coordinatize $N_{+} \cdot B_{-} \cong \n_{+}$ by $x_{\alpha}, \alpha \in \Delta^{+}(\g)$, with 
$x_{\alpha}$ dual to $E_{\alpha}$, and $LN_{+}$ by $x_{n,\alpha}, \, n \in \mathbb{Z}, 
\alpha \in \Delta^{+}(\g)$, with $x_{n,\alpha}$ dual to $E_{n,\alpha}$. Then $\sigma_{N_{+}}$ acts by sending
$x_{\alpha} \rightarrow x_{\sigg(\alpha)}$, $\partial_{x_{\alpha}} \rightarrow \partial_{x_{\sigg(\alpha)}}$. 

We introduce a Heisenberg algebra isomorphic to $\A$ with generators $\{ \tilde{a}_{n,\alpha}, \tilde{a}^{*}_{n,\alpha}, \mathbf{1} \}$ 
$,\alpha \in \Delta^{+}(\g),$ $n \in \mathbb{Z}$
and commutation relations 
\[
	[\tilde{a}_{n,\alpha},\tilde{a}_{m,\beta}] = [\tilde{a}^{*}_{n,\alpha},\tilde{a}^{*}_{m,\beta}] = 0
\]
\[
  	[\tilde{a}_{n,\alpha}, \tilde{a}^{*}_{m,\beta}] = \delta_{n,-m} \delta_{\alpha,\beta} \mathbf{1}
\]
This change of basis for $\A$ simply corresponds to the change of basis on $L\n_{+}$ from
a root basis to one on which $\sigLg$ acts diagonally. Under the identification $\tilde{a}_{n,\alpha} \rightarrow \partial x_{n,\alpha}, 
\tilde{a}^{*}_{n,\alpha} \rightarrow x_{-n,\alpha}$, $\A$ corresponds to the Weyl algebra of $LN_{+} \cdot LB_{-}$. 
In this basis for $\A$, $\sigM(\tilde{a}_{n,\alpha}) =  \tilde{a}_{n,\sigg(\alpha)}$, 
$\sigM(\tilde{a}^{*}_{n,\alpha}) = \tilde{a}^{*}_{n,\sigg(\alpha)}$.

Similarly, we introduce a Lie algebra isomorphic to $\H$ with generators $\{ \tilde{b}_{n,\alpha_{i}}, \mathbf{1} \}$
$n \in \mathbb{Z},$ $ i=1, \cdots, r$, and commutation relations
\[
	[\tilde{b}_{n,\alpha_{i}}, \tilde{b}_{m,\alpha_{j}}] = n (H_{\alpha_{i}}, H_{\alpha_{j}}) \delta_{n,-m} \mathbf{1}
\]
Here, $\sigpi(\tilde{b}_{n,\alpha_{i}}) = \tilde{b}_{n,\sigg(\alpha_{i})} $.

If $X \in \g$, let $R_{X}(x_{\alpha})$ be the polynomial vector field on $N_{+} \cdot B_{-}$ induced by $X$. We write
\[
	R_{X}(x_{\alpha}) = \sum_{\gamma \in \Delta^{+}(\g)} R^{\gamma}_{X}(x_{\alpha}) \partial x_{\gamma}
\]
We will use the notation $R^{M}_{X}$ to denote the vector of $\M$ given by the following expression
\[
	 \sum_{\gamma \in \Delta^{+}(\g)} R^{\gamma}_{X}(\tilde{a}^{*}_{0,\alpha}) \tilde{a}_{-1,\gamma} \vac
\]
where $R^{\gamma}_{X}(\tilde{a}^{*}_{0,\alpha})$ means that we replace $x_{\alpha}$ by $\tilde{a}^{*}_{0,\alpha}$ in the 
polynomial $R^{\gamma}_{X}$.
It follows from the untwisted Wakimoto construction (\cite{FF2}) that
\begin{align*}
	\iota(E_{-1,\alpha_{i}} \vac) &= R^{M}_{E_{\alpha_{i}}} \\
	\iota(H_{-1,\alpha_{i}} \vac) &= R^{M}_{H_{\alpha_{i}}} + \tilde{b}_{-1,\alpha_{i}} \\
	\iota(F_{-1,\alpha_{i}} \vac) &= R^{M}_{F_{\alpha_{i}}} + \lambda_{\alpha_{i}} \tilde{a}^{*}_{-1,\alpha_{i}} 
		- \tilde{a}^{*}_{0,\alpha_{i}} \tilde{b}_{-1,\alpha_{i}}
\end{align*}
where $\lambda_{\alpha}$ is some constant depending on the level $k$. It suffices to verify the lemma for these vectors, since
$\{ E_{\alpha_{i}}, F_{\alpha_{i}}, H_{\alpha_{i}}\} $ generate $\g$. 
Because the embedding $\g \rightarrow Vect(N_{+} \cdot B_{-})$ is equivariant, it follows that $\sigMpi(R^{M}_{X}) = R^{M}_{\sigg(X)}$. 
Therefore the embeddings of $\iota(E_{-1,\alpha_{i}} \vac), \iota(H_{-1,\alpha_{i}} \vac)$ are already equivariant, 
and it remains to prove that $\lambda_{\alpha_{i}} = \lambda_{\sigg(\alpha_{i})}$. 
To shorten notation, we denote the image of $X \in \Vg$ under $\iota$, $\iota(X)$, by $X^{\iota}$, and suppress the vacuum vector. Let 
\[
	S_{\alpha_{i}} = \sigMpi(F_{-1,\alpha_{i}}^{\iota}) - (\sigVg(F_{-1,\alpha_{i}}))^{\iota} = 
		(\lambda_{\alpha_{i}} - \lambda_{\sigg(\alpha_{i})}) \tilde{a}^{*}_{-1,\sigg(\alpha_{i})}
\]
Using the fact that $\iota$ is a vertex algebra homomorphism, and that $\Ymp(\sigMpi(v),z) = \sigMpi \Ymp(v,z) \sigMpi^{-1}$, we have:
\begin{align*}
	\Ymp & (\sigMpi(E_{-1,\alpha_{i}}^{\iota}),z) \Ymp(S_{\alpha_{i}},w) \\
		& \\ 
		\sim \, & \Ymp(\sigMpi(E_{-1,\alpha_{i}}^{\iota}),z) \Ymp ( \sigMpi(F_{-1,\alpha_{i}}^{\iota}),w ) \\
		& - \Ymp(\sigMpi(E_{-1,\alpha_{i}}^{\iota}),z) \Ymp((\sigVg(F_{-1,\alpha_{i}}))^{\iota} , w) \\
		& \\
		\sim \, & \sigMpi \Ymp(E_{-1,\alpha_{i}}^{\iota},z) , z) \Ymp(F_{-1,\alpha_{i}}^{\iota}, w) \sigMpi^{-1} \\
		& - \Ymp(E_{-1,\sigg(\alpha_{i})}^{\iota},z) \Ymp(F_{-1,\sigg(\alpha_{i})}^{\iota}, w ) \\
		& \\
		\sim \, & \frac{\sigMpi \Ymp(H_{-1,\alpha_{i}}^{\iota}, w) \sigMpi^{-1}}{(z-w)} + \frac{k(E_{\alpha_{i}}, F_{\alpha_{i}})}{(z-w)^{2}} \\
		& -  \frac{ \Ymp(H_{-1,\sigg(\alpha_{i})}^{\iota},w ) }{(z-w)}  -  \frac{k(E_{\sigg(\alpha_{i})}, F_{\sigg(\alpha_{i})})}{(z-w)^{2}} \\
		& \\
		\sim \, &  \frac{\Ymp(\sigMpi(H_{-1,\alpha_{i}}^{\iota}), w)}{(z-w)} -  \frac{ \Ymp(H_{-1,\sigg(\alpha_{i})}^{\iota},w ) }{(z-w)} \\
		\sim & 0 							
\end{align*}
However, we know from the Wakimoto construction that $R^{\alpha_{i}}_{E_{\alpha_{i}}} = -1$, so that $S_{\alpha_{i}} = 0$, which implies
that  $\lambda_{\alpha_{i}} = \lambda_{\sigg(\alpha_{i})} $ as desired. This proves that $\iota$ is equivariant. The equivariance
of $\kappa$ is shown similarly.

\end{proof}

\begin{theorem} 
$\fock$ has the structure of a $\ghats$ representation of level $k$ and highest weight $\chi$.
\end{theorem}

\begin{proof}

The homomorphism $\iota$ is $\sigma$ - equivariant, so that $\fock$ inherits a structure of $\sigma$ - twisted
$\Vg$ - module. Let $J^{1} \in \g_{i}, J^{2} \in \g$. Consider the twisted fields 
$\Ymps(\iota(J^{l}_{-1} \vac),z), \, l= 1,2$.
Using the associativity property for twisted modules (see \cite{BKT}) and the fact that 
$\iota$ is a vertex algebra homomorphism, we get (the vacuum $\vac$ is suppressed for brevity):
\begin{align*}
	 \Ymps & (\iota(J^{1}_{-1}),z)	\Ymps(\iota(J^{2}_{-1} ),w)  \\
	& \\
		 \sim & \,  \frac{ z^{-i/N} w^{i/N} \Ymps(\iota(J^{1}_{-1} )[0] \iota(J^{2}_{-1}),w )}{(z-w)} \\
		 & + \Ymps(\iota(J^{1}_{-1})[1] \iota(J^{2}_{-1} ),w ) \partial_{w} \frac{z^{-i/N}w^{i/N}}{(z-w)} \\
	& \\
  		\sim & \, \frac{ z^{-i/N} w^{i/N} \Ymps(\iota(J^{1}_{0} J^{2}_{-1}) ,w )}{(z-w)} \\
		 & + \Ymps(\iota(J^{1}_{1} J^{2}_{-1} ),w ) \partial_{w} \frac{z^{-i/N}w^{i/N}}{(z-w)} \\
	& \\
		\sim & \, \frac{ z^{-i/N} w^{i/N} \Ymps(\iota([J^{1}, J^{2}]_{-1}) ,w )}{(z-w)} \\
		&  + k(J_{1},J_{2}) \partial_{w} \frac{z^{-i/N}w^{i/N}}{(z-w)} \\
\end{align*}
This implies that the fields $\Ymps(\iota(J^{l}_{-1}),z)$ have the same OPE as ( \ref{twistedOPE} ), and hence 
their Fourier coefficients generate $\ghats$.
\end{proof}

We will henceforth use the notation $\waks$ for $\fock$. We call $\waks$ the twisted Wakimoto module of 
level $k$ and highest weight $\chi$. 

\begin{theorem}	\label{CharThm}
The $\ghats$ module $\waks$ has the same formal character as $\verm$ - the Verma module of level $k$ with highest weight
$\chi$.  
\end{theorem}

\begin{proof}

The Cartan subalgebra $\h^{\sigma} = \h_{0} \oplus \mathbb{C} \cdot K \oplus \mathbb{C} \cdot d$ 
of $\ghats$ is spanned by $\{ H_{0,(0,a)}, K, d \}, \, a = 1, \cdots , dim(\h_{0})$. When the level $k$ is understood, we write elements 
$\gamma \in \h^{\sigma,*}$ in the form
$\bar{\gamma} + n \delta, \, n \in \frac{1}{N} \mathbb{Z}$, where $\bar{\gamma}$ is the orthogonal projection of $\gamma$ to $\h^{*}_{0}$, 
and $\delta$ is dual to $d$.

Unraveling the homomorphism $\iota$, we see that $H_{0,(0,a)}$ acts on $\fock$
as the operator 
\[
	\sum_{(j,\alpha) \atop \alpha \in \Delta^{+}_{j}} \sum_{n \in - \frac{j}{N} + \mathbb{Z}} 
	   \alpha(H_{0,a}) :a^{*}_{n,(j,\alpha)}a_{-n,(j,\alpha)}: + b_{0,(0,a)}
\]
$\fock$ has a basis of monomials of the form 
\[
		 \prod_{(j,\alpha), \alpha \in \Delta^{+}_{j} \atop n \in - \frac{j}{N} + \mathbb{Z}, n \leq 0} 
			(a^{*}_{n,(j,\alpha)})^{\lambda_{n,(j,\alpha)}}
		\prod_{(j,\alpha), \alpha \in \Delta^{+}_{j} \atop n \in \frac{j}{N} + \mathbb{Z}, n < 0 } 
			(a_{n,(j,\alpha)})^{\mu_{n,(j,\alpha)}}	
		 \prod_{(i,a) \atop n \in \frac{i}{N} + \mathbb{Z}, n < 0  } (b_{n,(i,a)})^{\tau_{n,(i,a)}}
\]
where $\lambda_{n,(j,\alpha)}$ and $\tau_{n,(i,a)}$ are non-negative integers. 
Now
\begin{align*}
	[H_{0,(0,a)}, a^{*}_{n,(j,\alpha)}] &= - \alpha(H_{0,(0,a)}) \\
	[H_{0,(0,a)}, a_{n,(j,\alpha)}] &= \alpha(H_{0,(0,a)}) \\
	[H_{0,(0,a)}, b_{n,(i,c)}] &= 0 
\end{align*}
An easy computation shows that the above monomial lies in the
\[	\chi +
	\sum_{(j,\alpha), \alpha \in \Delta^{+}_{j} \atop n \in - \frac{j}{N} + \mathbb{Z}, n \leq 0} \lambda_{n,(j,\alpha)} (- \alpha + n \delta) + 
	\sum_{(j,\alpha), \alpha \in \Delta^{+}_{j} \atop n \in \frac{j}{N} + \mathbb{Z}, n < 0} \mu_{n,(j,\alpha)} (\alpha + n \delta) + 
	\sum_{(i,a) \atop n \in \frac{i}{N} + \mathbb{Z}, n < 0} n \tau_{n,(i,a)} \delta 
\]
weight space. 
The formal character of $\fock$ is therefore given by:
\begin{align}
	& Ch(\fock) =  \nonumber \label{char} \\
	& e^{\chi} \prod_{(j,\alpha), \alpha \in \Delta^{+}_{j} \atop n \in - \frac{j}{N} + \mathbb{Z}, n \geq 0} 
				\frac{1}{1-e^{-(\alpha + n \delta)}} 
		    \prod_{(j,\alpha), \alpha \in \Delta^{+}_{j} \atop n \in \frac{j}{N} + \mathbb{Z}, n > 0} 
				\frac{1}{1-e^{-(-\alpha + n \delta)}}
		    \prod_{(i,a) \atop n \in \frac{i}{N} + \mathbb{Z}, n > 0} \frac{1}{1-e^{- n \delta}}  \\
	& = e^{\chi} \prod_{\beta \in \Delta_{+}(\ghats)} \frac{1}{(1-e^{- \beta})^{mult(\beta)}} \nonumber
\end{align}
which is the character of $\verm$ - the Verma module.
\end{proof}

\section{Critical Level}

\subsection{Restricted Wakimoto Modules}	\label{ResWak}

The $\Hs$ module $\pi^{0,\sigma}_{\chi}$ has a proper submodule $\bar{\pi}^{0,\sigma}_{\chi}$ spanned by all non-constant monomials
in $b_{n,(i,a)}, \, n < 0 $, and in the notation of section \ref{FreeBoson} 
$\pi^{0,\sigma}_{\chi} / \bar{\pi}^{0,\sigma}_{\chi} \cong \mathbf{C}_{0,\chi} $, where we have extended $\mathbf{C}_{0,\chi}$
by $0$ from $\Hs_{+} \oplus \mathbb{C} \cdot \mathbf{1}$ 

Thus when $k = - \dcox$, the fields $\Ymps$ can be made to act on the smaller Fock-space
\[
	\rwaks = \Ms \otimes \mathbf{C}_{0, \chi}
\]
We obtain

\begin{theorem} 
	$\rwaks$ has the structure of a $\ghats$ representation of highest weight $\chi$ and level
	$- \dcox$
\end{theorem}

We will henceforth refer to the module $\rwaks$ as the restricted Wakimoto module of critical level, following
\cite{FF2}.

\subsection{Restricted Contragredient Wakimoto Modules}

We will also make use of the contragredient Wakimoto module $\crwaks$. In this section
we give an explicit realization of $\crwaks$ in terms of a Fock space, analogous to that of $\rwaks$. 
Let $\omega$ denote the Chevalley anti-involution of $\ghats$. $\omega$ fixes $\h_{0} \oplus \mathbb{C} \cdot d \oplus \mathbb{C} \cdot K$
, and maps  $E_{n,(i,\alpha)}$ to a multiple of $E_{-n,(-i,-\alpha)}$. 

Given an $\h$ - diagonalizable category $\mathcal{O}$ module $V$ with finite-dimensional weight spaces, we define 
its contragredient module $V^{*}$ to be the restricted dual of $V$ on which $\ghats$ acts by
\[
	x \cdot \phi (v) = \phi (\omega(x) \cdot v)
\]
where $x \in \ghats, v \in V, \phi \in V^{*}$.

Let $\overline{\As}_{+} \subset \As$ denote the abelian subalgebra generated by $a_{n,(j,\alpha)}, a^{*}_{m,(k,\beta)}$, $ \, n \leq 0, m < 0$.
Let $\Msbar = Ind^{\As}_{\overline{\As}_{+}} \mathbf{C}$, where $\mathbf{C}$ denotes the trivial representation of $\overline{\As}_{+}$.
$\Msbar$ has a basis of monomials of the form 
\[
	Z(\lambda,\mu)=
	\prod_{n_{r},(j_{r},\alpha_{r}) \atop n_{r} \leq 0} (a_{n_{r},(j_{r},\alpha_{r})})^{\lambda_{n_{r},(j_{r}, \alpha_{r})}} 
	\prod_{n_{s},(j_{s},\alpha_{s}) \atop n_{s} < 0} (a^{*}_{n_{s},(j_{s},\alpha_{s})})^{\mu_{n_{s},(j_{s}, \alpha_{s})}} \vac
\]
We define $\overline{Z(\lambda,\mu)} \in \Ms$ by
\[
	\overline{Z(\lambda,\mu)}=
	\prod_{n_{r},(j_{r},\alpha_{r}) \atop n_{r} \leq 0} (a^{*}_{n_{r},(j_{r},\alpha_{r})})^{\lambda_{n_{r},(j_{r}, \alpha_{r})}} 
	\prod_{n_{s},(j_{s},\alpha_{s}) \atop n_{s} < 0} (a_{n_{s},(j_{s},\alpha_{s})})^{\mu_{n_{s},(j_{s}, \alpha_{s})}} \vac
\]
We now define a non-degenerate bilinear pairing $\langle , \rangle: \Msbar \otimes \Ms \rightarrow \mathbb{C}$ by
\[
	\langle \overline{Z(\lambda,\mu)}, Z(\lambda,\mu) \rangle = 
		\prod_{n_{r},(j_{r},\alpha_{r}) \atop n_{r} \leq 0} \lambda_{n_{r},(j_{r}, \alpha_{r})} !
		\prod_{n_{s},(j_{s},\alpha_{s}) \atop n_{s} < 0} \mu_{n_{s},(j_{s}, \alpha_{s})} !
\]
For $v \in \Msbar, w \in \Ms$, we have
\begin{align*}
	\langle v, a_{n,(j,\alpha)} \cdot w \rangle & = \langle a_{-n,(-j,\alpha)} v, w \rangle \\
	\langle v, a^{*}_{n,(j,\alpha)} \cdot w \rangle & = \langle -a^{*}_{-n,(-j,\alpha)} \cdot v, w \rangle
\end{align*} 
Let $T: End(\Ms) \rightarrow End(\Msbar)$ denote the ``transpose'' map defined by
\begin{align*}
	T(a_{n,(j,\alpha)}) & = a_{-n,(-j,\alpha)} \\
	T(a^{*}_{n,(j,\alpha)}) & = -a^{*}_{-n,(-j,\alpha)}
\end{align*}
and the property $T(A \cdot B) = T(B) \cdot T(A), \, A,B \in End(\Ms)$. We have 
\begin{equation}
	\langle v, A \cdot w \rangle = \langle T(A) \cdot v, w \rangle	\label{adjointrelation}
\end{equation}
It follows from (\ref{adjointrelation}) that if $x \in \ghats$ then 
\[
	\langle v,\omega(x) \cdot w \rangle = \langle T(\omega(x)) \cdot v, w \rangle 
\]
Thus, if we identify $\Msbar$ with the restricted dual $\M^{\sigma,*}$ via the pairing $\langle , \rangle$, then
the action 
\[
	x(v) = T(\omega(x)) \cdot v
\]
identifies $\Msbar$ with the $\ghats$ module contragredient to $M = \overline{W}^{\sigma}_{0}$. A simple computation shows that
at the level of fields, the formula for $E_{j,\alpha}(z)$ acting on $\Msbar$ is the old formula for $\omega(E_{j,\alpha})(z)$
acting on $\Ms$, possibly with some sign changes among the monomials. Let  
\[
	\crwaks = \Msbar \otimes \mathbb{C}_{0,\chi}
\]
In the same manner as in section \ref{ResWak}, we extend the $\ghats$ action on $\Msbar$ to 
$\crwaks$. 

\section{Proof of the Kac-Kazhdan Conjecture in the Twisted Case}

In this section, let $(,)$ denote the inner product on $\h^{\sigma, *}$ defined as in (\cite{K1}).
At the critical level $k = - \dcox$, a $\ghats$ highest weight $\chi \in \h^{\sigma,*}$ automatically satisfies 
the Kac-Kazhdan equation (see \cite{KK})
\begin{equation}
	2( \chi + \rho, \beta) = n (\beta,\beta) \label{KKeqn}
\end{equation}
corresponding to imaginary roots $\beta \in \Delta^{+}_{im}(\ghats)$. This implies that the Verma module $M_{\chi}$ is reducible. 
Let $L_{\chi}$ denote the irreducible quotient of $M_{\chi}$. In \cite{KK}, it is conjectured
that for $\chi$ generic, i.e. such $\chi$ that do not satisfy the equation (\ref{KKeqn}) for
any real roots $\beta$, the formal character of $L_{\chi}$ is
\[
	Ch(L_{\chi}) =  e^{\chi} \prod_{\alpha \in \Delta^{re}_{+}(\ghats)} \frac{1}{(1-e^{- \alpha})}
\]
In the case of untwisted affine algebras, the Kac-Kazhdan conjecture was proved in \cite{F1}. In this
section, we treat the twisted case using the same approach. The conjecture follows immediately from the following
theorem:

\begin{theorem}	\label{IrredThm}
	The $\ghats$ representation $\rwaks$ is irreducible for generic critical $\chi$. 
\end{theorem}

\subsection{Proof of Theorem \ref{IrredThm}}

The proof of the theorem will require several steps. \\

There exist operators  $\bar{E}_{n,(j,\alpha)}, n \in \frac{j}{N} + \mathbb{Z}, \alpha \in \Delta^{+}_{j}$, acting on $\rwaks$,
which generate a Lie algebra isomorphic to $\Lns$ and commute with $E_{n,(j,\alpha)},  n \in \frac{j}{N} + \mathbb{Z}, \alpha \in 
\Delta^{+}_{j}$. They arise from the right action of 
$\Lgs$ on the big cell of the twisted semi-infinite Flag Manifold $LG^{\sigma}/ LB^{\sigma}_{-}$. 
Let 
\[
	\bar{E}_{j,\alpha}(z) = \sum_{n \in \frac{j}{N} +\mathbb{Z}} \bar{E}_{n,(j,\alpha)} z^{-n-1}
\]
The nilpotence of $\Lns$ implies that the $\bar{E}_{n,(j,\alpha)}$'s are given by formulas of the form
\begin{equation}
	\bar{E}_{j,\alpha}(z) = a_{j,\alpha}(z) + \sum_{(k,\beta) \atop 0 \leq k < N, \beta \in \Delta^{+}_{k}} a_{k,\beta}(z) 
			P_{(j,\alpha)}^{(k,\beta)}(a^{*}_{(l,\gamma)}(z)) 	\label{ebarformula}
\end{equation}
where $\beta \succ \alpha$ and $\beta \succ \gamma$ in the root lattice of $\g_{0}$, and $P_{(j,\alpha)}^{(k,\beta)}$ are polynomials.

The operators $\bar{E}_{n,(j,\alpha)}, n < 0$ generate a Lie algebra isomorphic to  
\[
L_{-} \n^{\sigma}_{+}= \Lns \cap \g \otimes t^{- \frac{1}{N}} \mathbb{C}[[t^{- \frac{1}{N}}]]
\]
Let $F$ be the subspace of $\rwaks$ generated by $\bar{E}_{n,(j,\alpha)}, n < 0$  from the vacuum. We have a map
\[
	\jmath: F \otimes \mathbb{C}[a^{*}_{n,(j,\alpha)}]_{n \leq 0} \rightarrow \rwaks
\]
sending a vector of the form 
\[
	\bar{E}_{n_{1},(j_{1},\alpha_{1})} \cdots \bar{E}_{n_{k},(j_{k},\alpha_{r})} \vac \otimes a^{*}_{m_{1},(k_{1},\beta_{1})} \cdots 
	a^{*}_{m_{s},(k_{s},\beta_{s})} \vac
\]
to the element of $\rwaks$ given by
\[ 
	\bar{E}_{n_{1},(j_{1},\alpha_{1})} \cdots \bar{E}_{n_{k},(j_{k},\alpha_{r})}
	a^{*}_{m_{1},(k_{1},\beta_{1})} \cdots 
	a^{*}_{m_{s},(k_{s},\beta_{s})} \vac 
\]

\begin{lemma}
	The map $\jmath$ is an isomorphism
\end{lemma}
\begin{proof}
It is clear that $\jmath$ is injective.
The expression (\ref{ebarformula}) for $\bar{E}_{j,\alpha}(z)$ implies that we can solve for $a_{j,\alpha}(z)$ in terms
of $a^{*}_{k,\beta}(z)$ and $\bar{E}_{i,\gamma}(z)$. In other words, 
\[
	a_{j,\alpha}(z) = \bar{E}_{j,\alpha}(z) + \sum  \bar{E}_{k,\beta}(z) R^{(k,\beta)}_{(j,\alpha)}(a^{*}_{i,\gamma}(z))
\] 
where in the summation, $\beta \succ \alpha$ in the root lattice of $\g_{0}$, and $R^{(k,\beta)}_{(j,\alpha)}$
are polynomials. Now, $\rwaks$ has a basis of monomials of
the form 
\[
	a_{n_{1},(j_{1},\alpha_{1})} \cdots a_{n_{r},(j_{r},\alpha_{s})}  
		a^{*}_{m_{1},(k_{1},\beta_{1})} \cdots a^{*}_{m_{s},(k_{s},\beta_{s})} \vac
\]
where $n_{p} < 0, m_{q} \leq 0$. Thus we see that $\jmath$ is surjective.

\end{proof}

We use the map $\jmath$ to identify $F \otimes \mathbb{C}[a^{*}_{n,(j,\alpha)}]_{n \leq 0}$ and $\rwaks$. The 
operators $\bar{E}_{n,(j,\alpha)}$ commute with the action of $\Lns \in \ghats$, and so under this identification,
if $P \in F, Q \in \mathbb{C}[a^{*}_{n,(j,\alpha)}]_{n \leq 0}$ 
\[
	E_{m,(k,\beta)} \cdot (P \otimes Q) = P \otimes (E_{m,(k,\beta)} \cdot Q)
\]
Let
\[
 	L_{+} \n^{\sigma}_{+} = \Lns \cap  \g \otimes t^{\frac{1}{N}} \mathbb{C}[[t^{\frac{1}{N}}]] \subset \ghats
\]
Then
$L_{+} \n^{\sigma}_{+}$ acts co-freely on $\mathbb{C}[a^{*}_{n,(j,\alpha)}]_{n \leq 0}$, i.e. 
given $X \in \mathcal{U} (L_{+} \n^{\sigma}_{+})$ (here $\mathcal{U}()$ stands for enveloping algebra), 
there is a unique $Q \in \mathbb{C}[a^{*}_{n,(j,\alpha)}]_{n \leq 0}$
such that $X \cdot Q = \vac $. Now suppose that $v \in \rwaks$ is a singular vector. Then in
particular, it must be killed by $L_{+} \n^{\sigma}_{+}$, and so must be of the form 
$P \otimes \vac, P \in F$. Thus all potential singular vectors must lie in the subspace
$F \otimes \vac \subset F \otimes \mathbb{C}[a^{*}_{n,(j,\alpha)}]_{n \leq 0}$. \\

\noindent The module $\rwaks$ is graded by $\Lambda(\ghats)$ (if $v \in \rwaks$ belongs to the
$\mu$ weight space, $deg(v) = \chi - \mu$). Using the notation of Theorem \ref{CharThm}, we write
$deg(v)=\beta = \bar{\beta} + m \delta$, where $\beta \in \Lambda(\g), m \in \frac{1}{N} \mathbb{Z}$.

\begin{lemma} 	\label{DegLemma}
 	If $v \in F \otimes \vac$, $deg(v)= \beta$, then $\bar{\beta} \ne 0$
\end{lemma}
\begin{proof}
	The degree of the operator $\bar{E}_{n,(j,\alpha)}, n < 0, \alpha \in \Delta^{+}(j)$ is 
$\alpha + n \delta$. It follows that repeated 
application of such operators to the vacuum can never result in a vector with degree $m \delta$. 
\end{proof}

We can apply the same line of argument to the contragredient module $\crwaks$. The operators $E_{n,(j,-\alpha)}$
$n > 0, \alpha \in \Delta^{+}_{j}$ generate the Lie algebra 
\[
	L_{+}n^{\sigma}_{-} = \n_{-} \otimes t^{\frac{1}{N}} \mathbb{C}[[t^{\frac{1}{N}}]] \subset \ghats
\]
The operators $\bar{E}_{n,(j,\alpha)}$ act on $\crwaks$ as before, and commute with the action of $L_{+} \n^{\sigma}_{-}$.
Let $\bar{F}$ denote the subspace generated by $\bar{E}_{n,(j,\alpha)}, \, n \leq 0$ from the vacuum in $\crwaks$. We can show
as before that $\crwaks = \bar{F} \otimes \mathbb{C}[a^{*}_{n,(j,\alpha)}]_{n < 0}$, that if $E_{n,(j,\alpha)} \in L_{+} \n^{\sigma}_{-},
P \in \bar{F}, Q \in \mathbb{C}[a^{*}_{n,(j,\alpha)}]_{n < 0} $, then
\[
	E_{n,(j,\alpha)} \cdot (P \otimes Q) = P \otimes (E_{n,(j,\alpha)} \cdot Q)
\]
and that the action of $L_{+} \n^{\sigma}_{-}$ on $ \mathbb{C}[a^{*}_{n,(j,\alpha)}]_{n < 0} $ is co-free. Thus, any singular vector
$v \in \crwaks$ lies in $\bar{F} \otimes \vac$, and can be shown to have degree $deg(v) = \beta, \, \bar{\beta} \ne 0$. \\

We are now ready to prove Theorem \ref{IrredThm}

\begin{proof}
The formal character $Ch(V)$ of a category $\mathcal{O}$ $\ghats$ module $V$ can be uniquely expressed as 
a linear combination of characters of irreducible highest-weight modules $Ch(L_{\psi}), \psi \in \h^{\sigma,*}$. 
Let $C_{V}^{\psi} \in \mathbb{Z}$ denote the coefficient with which $Ch(L_{\psi})$ occurs in $Ch(V)$.

Suppose that $\rwaks$ has a proper submodule. Then either $\rwaks$ or $\crwaks$ has a singular vector $v \ne \vac$. Suppose
first that $v \in \rwaks$, and that $v$ has weight $\mu$. It follows that $C_{\rwaks}^{\mu} > 0$. Now, 
$\rwaks$ is a quotient of $W^{\sigma}_{- \dcox, \chi}= \Ms \otimes \pi^{0,\sigma}_{\chi}$, which by Theorem \ref{CharThm} 
has the same character as the Verma module $M_{-\dcox, \chi}$. It follows that $C_{M_{-\dcox,\chi}}^{\mu} > 0$. 
By Lemma \ref{DegLemma}, $\chi - \mu = \beta, \, \bar{\beta} \ne 0$. Proposition (4.1) from \cite{KK} then
implies that $M_{-\dcox, \chi}$ contains a singular vector of weight $\mu$, and hence degree $\bar{\beta}+m \delta$, contradicting
the following Lemma. \\

\noindent The case where $v \in \crwaks$ is treated in the same way.
\end{proof}

\begin{lemma}
Let $\chi \in \h^{\sigma,*}$ be a generic critical highest weight, and $v \in M_{-\dcox,\chi}$ a singular vector. Then
$deg(v) = m \delta , m \in \frac{1}{N} \mathbb{Z}$
\end{lemma}
\begin{proof}
Suppose that $v$ has weight $\mu$. By Theorem 2 in \cite{KK}, there exists a sequence of 
positive roots $\gamma_{1}, \cdots, \gamma_{k} \in \Delta^{+}(\ghats)$, and a sequence of positive integers
$n_{i}, \cdots, n_{k}$ such that 
\[
	2(\chi + \rho -n_{1}\gamma_{1} - \cdots - n_{i-1} \gamma_{i-1},\gamma_{i}) = n_{i}(\gamma_{i},\gamma_{i})
\]
for $i = 1, \cdots, k$ and $\chi - \mu = \sum_{i=1}^{k} n_{i}\gamma_{i}$. 
Now $\chi - \mu = \bar{\beta} + m\delta , \, \bar{\beta} \ne 0$, so $\gamma_{r} \in \Delta^{+}_{re}(\ghats)$ for some
$r$. Let us choose the smallest such $r$. We have that
\[
	2(\chi + \rho - n_{1}\gamma_{1} - \cdots - n_{r-1} \gamma_{r-1}, \gamma_{r}) = n_{r}(\gamma_{r},\gamma_{r})
\] 
and since real and imaginary roots are orthogonal, $2(\chi + \rho,\gamma_{r}) = n_{r} (\gamma_{r}, \gamma_{r})$,
contradicting the assumption that $\chi$ is generic. 
\end{proof}

\section{Wakimoto Realization of $A_{2}^{(2)}$}

In this section we give explicit realizations of the twisted affine algebra $A^{(2)}_{2}$.
Let $R_{i,j}$ denote the $3 \times 3$ matrix with a $1$ in position $(i,j)$.  

Here $\g= \mathfrak{sl_{3}}$, where we take the standard realization as traceless $3 x 3$
matrices. Under $\sigma$, $\g = \g_{0} \oplus \g_{1}$, where
$\g_{0} = \mathfrak{sl_{2}}$. As a $\g_{0}$ - representation, $\g_{1}$ is isomorphic
to the irreducible $5$ - dimensional representation. Let $\alpha$ denote the highest
root of $\mathfrak{sl_{2}}$. We choose a basis as follows
\begin{align*}
E_{0,\alpha} &= R_{(1,1)}+R_{(2,3)} \\ 
E_{1,\alpha} &= R_{(1,1)}-R_{(2,3)} \\
E_{1,2 \alpha} &=-2 R_{(1,3)} \\  
H_{0,1} &= R_{(1,1)}-R_{(3,3)} \\  
H_{1,1} &= R_{(1,1)}-2R_{(2,2)}+R_{(3,3)} \\  
E_{0,-\alpha} &= R_{(2,1)}+R_{(3,2)} \\ 
E_{1,-\alpha} &= R_{(2,1)}-R_{(3,2)} \\  
E_{1,-2\alpha} &= -2 R_{(3,1)}
\end{align*}

The Wakimoto realization in this case looks as follows (we omit the fields $E_{1,2 \alpha}(z), E_{1,-2 \alpha}(z)$ as
these can be calculated via the OPE from the other generators ):

\begin{align*}
	E_{(0,\alpha)}(z) &= -a_{(0,\alpha)}(z) -1/2 :a^{*}_{(1,\alpha)}(z) a_{(1,2 \alpha)}(z): \\
	E_{(1,\alpha)}(z) &= -a_{(1,\alpha)}(z) +1/2 :a^{*}_{(0,\alpha)}(z) a_{(1,2 \alpha)}(z): \\
	H_{(0,1)}(z) &= - :a^{*}_{(0,\alpha)}(z) a_{(0,\alpha)}(z): - :a^{*}_{(1,\alpha)}(z) a_{(1,\alpha)}(z):
		-2 :a^{*}_{(1,2 \alpha)}(z) a_{(1,2 \alpha)}(z): + b_{(0,1)}(z) \\
	H_{(1,1)}(z) &= -3 :a^{*}_{(0,\alpha)}(z) a_{(0, \alpha)}(z): -3 :a^{*}_{(1,\alpha)}(z) a_{(0, \alpha)}(z): + b_{(1,1)}(z) \\
	E_{(0,- \alpha)}(z) &= 1/2 :(a^{*}_{(0,\alpha)}(z))^{2} a_{(0, \alpha)}(z): + 
		3/2:(a^{*}_{(1,\alpha)}(z))^{2} a_{(0, \alpha)}(z): \\
		& +2 :a^{*}_{(0,\alpha)}(z) a^{*}_{(1,\alpha)}(z) a_{(1, \alpha)}(z): 
		- 2 :a^{*}_{(1,2 \alpha)}(z) a_{(1, \alpha)}(z): + \\ 
		& 1/4 :a^{*}_{(1,\alpha)}(z) (a^{*}_{(0,\alpha)}(z))^{2} a_{(1, 2 \alpha)}(z): 
		-1/4 :(a^{*}_{(1,\alpha)}(z))^{3} a_{(1, 2 \alpha)}(z): + \\
		& : a^{*}_{(0,\alpha)}(z) a^{*}_{(1,2 \alpha)}(z) a_{(1, 2 \alpha)}(z): + 
			(-1 - 2k) \partial_{z} a^{*}_{0,\alpha}(z) - a^{*}_{0,\alpha}(z)b_{(0,1)}(z) \\
	E_{(1,- \alpha)}(z) &= 2 :a^{*}_{(1,2 \alpha)}(z)a_{(0, \alpha)}(z): + 
			2:a^{*}_{(0,\alpha)}(z)a^{*}_{(1,\alpha)}(z)a_{(0, \alpha)}(z): \\
			    & + 3/2 :(a^{*}_{(0,\alpha)}(z))^{2} a_{(1, \alpha)}(z): + 
				1/2 :(a^{*}_{(1,\alpha)}(z))^{2} a_{(1, \alpha)}(z):  \\
			    & 	+1/4 :(a^{*}_{(0,\alpha)}(z))^{3}a_{(1, 2 \alpha)}(z): 
				 -1/4 :a^{*}_{(0,\alpha)}(z)(a^{*}_{(1,\alpha)}(z))^{2}a_{(1, 2 \alpha)}(z): \\
			    &	+ :a^{*}_{(1,\alpha)}(z)a^{*}_{(1,2 \alpha)}(z) a_{(1, 2 \alpha)}(z):
				+ (- 1 - 2k) \partial_{z} a^{*}_{1,\alpha}(z) - a^{*}_{1,\alpha}(z)b_{(1,1)}(z) \\
\end{align*}

\textbf{Remark:} After this work was completed, we were informed by Prof. T. Takebe that some of our results had been obtained
independently by Prof. G. Kuroki

\end{document}